\documentclass[12pt,reqno]{amsart}
\usepackage[latin1]{inputenc}
\usepackage{amsmath}
\usepackage{amssymb}
\usepackage{amscd}
\usepackage{natbib}
\usepackage{mathrsfs}
\usepackage[pdftex]{graphicx}
\usepackage{graphicx}
\usepackage[OT2, T1]{fontenc}
\usepackage[russian,english]{babel}
\usepackage[babel=true]{csquotes}
\usepackage{chngcntr}
\counterwithout{figure}{section}
\counterwithout{figure}{subsection}

\newcommand{\ix}[0]{\ensuremath{\mathbf{x}}}

\newcommand{\ka}[0]{\ensuremath{\mathbf{k}}}

\newcommand{\zero}[0]{\ensuremath{\mathbf{0}}}

\newcommand{\R}[0]{\ensuremath{\mathbb{R}}}
\newcommand{\T}[0]{\ensuremath{\mathbb{T}}}
\newcommand{\Z}[0]{\ensuremath{\mathbb{Z}}}

\DeclareMathOperator{\esssup}{ess\,sup} 

\newtheorem{theo}{Theorem}[section]{\bf}{\it}
\newtheorem{lemm}[theo]{Lemma}{\bf}{\it}
\newtheorem{defi}[theo]{Definition}{\bf}{\it}
\newtheorem{cor}[theo]{Corollary}{\bf}{\it}
\newtheorem{rmq}[theo]{Remark}{\bf}{\it}

\author{Alexandre Boritchev}

\address{University of Lyon
\\
CNRS UMR 5208
\\
University Claude Bernard Lyon 1
\\
Institut Camille Jordan
\\
43 Blvd. du 11 novembre 1918
\\
69622
\\
VILLEURBANNE CEDEX
\\
FRANCE
\\
E-mail:
\\
alexandre.boritchev@gmail.com
}
\title[Decaying subcritical Burgulence]{Decaying turbulence for the fractional subcritical Burgers equation}
\date{\today}
\begin{document}


\maketitle

\bigskip
\textbf{Abstract.}\ We consider the fractional unforced Burgers equation in the one-dimensional space-periodic setting:
\begin{equation} \nonumber
\frac{\partial u}{\partial t}+(f(u))_x +\nu \Lambda^{\alpha} u= 0,\quad t \geq 0,\ \ix \in \T^d=(\R/\Z)^d.
\end{equation}
Here $f$ is strongly convex and satisfies an additional growth condition, $\Lambda=\sqrt{-\Delta}$, $\nu$ is small and positive, while $\alpha \in (1,\ 2)$ is a constant in the \textit{subcritical} range.
\\ \indent
For solutions $u$ of this equation, we generalise the results obtained for the case $\alpha=2$ (i.e. when $-\Lambda^{\alpha}$ is the Laplacian) in \cite{BorD}. We obtain sharp estimates for the time-averaged Sobolev norms of $u$ as a function of $\nu$. These results yield sharp estimates for natural analogues of quantities characterising the hydrodynamical turbulence, namely the averages of the increments and of the energy spectrum. In the inertial range, these quantities behave as a power of the norm of the relevant parameter, which is respectively the separation $\ell$ in the physical space and the wavenumber $\ka$ in the Fourier space.
\\ \indent
The form of all estimates is the same as in the case $\alpha=2$; the only thing that changes (except implicit constants) is that $\nu$ is replaced by $\nu^{1/(\alpha-1)}$.

\tableofcontents

\tableofcontents

\section{Introduction} \label{intro}

\subsection{Burgers turbulence} \label{burg}
\smallskip
\indent
The Burgers equation
\begin{equation} \label{classBurgers}
\frac{\partial u}{\partial t} + u \frac{\partial u}{\partial x} - \nu \frac{\partial^2 u}{\partial x^2} = 0,
\end{equation}
where $\nu>0$ is a constant, and its multidimensional generalisations, are very popular physical models: see the review \cite{BK07} and references therein. On a formal level, this equation looks like a pressureless one-dimensional model for the incompressible Navier-Stokes equations \cite{Pol95}. In the turbulent regime, i.e. for $\nu \ll 1$, the solutions of the Burgers equation display non-trivial small-scale intermittent behaviour, called decaying Burgers turbulence or \enquote{Burgulence} \cite{Bur74,Cho75,Kid79}
\\ \indent
To fix the ideas, let us now consider the space-periodic setting, i.e. $x \in S^1=\R/\Z$, and an initial condition of order $1$ and recall some standard qualitative facts about the behaviour of the solutions of (\ref{classBurgers}); for details, see \cite[Section 1]{BorD}. After a time of order $1$, these solutions have an $N$-wave behaviour in the limit $\nu \rightarrow 0$. In other words, $u(t,\cdot)$ displays negative jump discontinuities separated by smooth regions where the derivative is positive \cite{Eva08}. For $0<\nu \ll 1$ the solutions are still highly intermittent: the jump discontinuities become layers of width of order $\nu$ where the derivative is negative of order $\nu^{-1}$. These layers are called {\it cliffs} \cite{Fri95}.
\\ \indent
On a physical level of rigour, the arguments given above imply two results for the small-scale behaviour of the solutions:
\begin{itemize}
\item On one hand, for $\nu$ small and for $1 \ll k \ll \nu^{-1}$, the energy-type quantities $\frac{1}{2} |\hat{u}(k)|^2$ (where $\hat{u}(k)$ is the $k$-th Fourier coefficient) behave, in average, as $k^{-2}$ \cite{Cho75,FouFri83,Kid79,Kra68}.
\item On the other hand, the structure functions
$$
S_p(\ell)=\int_{S^1}{|u(x+\ell)-u(x)|^p\ dx}
$$
behave as $\ell^{\max(1,p)}$ for $\nu \ll \ell \ll 1$: in other words, we have a \textit{bifractal behaviour}: see \cite{AFLV92} and \cite[Chapter 8]{Fri95}.
\end{itemize}
\indent
In the description above, we see that the length scale of the system is of order $\nu$. Heuristically, this can be justified by looking at the form of the equation. First we assume that the solution $u$ is of order $1$ and we ignore the term $u_t$. We denote the length scale by $\ell$ and therefore taking the derivative $\partial_x$ amounts to multiplying by $\ell^{-1}$. Then we obtain that 
$uu_x \sim \ell^{-1}$ and $\nu u_{xx} \sim \nu \ell^{-2}$, which yields that $\ell \sim \nu$.

\subsection{Fractional Burgers equation} \label{introturb}
\smallskip
\indent
Now we consider the fractional Burgers equation
\begin{equation} \label{fracBurgers}
\frac{\partial u}{\partial t} + u \frac{\partial u}{\partial x} + \nu \Lambda^{\alpha} u = 0, \quad x \in S^1=\R/\Z,\ \alpha \geq 0,
\end{equation}
where $\nu$ is a constant in $(0,1]$, $\Lambda=\sqrt{-\Delta}$ and $f$ is $C^{\infty}$-smooth and strongly convex, i.e. $f$ satisfies the property
\begin{equation} \label{strconvex}
f''(y) \geq \sigma > 0,\quad y \in \R.
\end{equation}
For the sake of simplicity, we only consider solutions to (\ref{fracBurgers})-(\ref{strconvex}) with zero space average for fixed $t$:
\begin{equation} \label{zero}
\int_{S^1}{u(t,x) dx}=0,\quad \forall t \geq 0.
\end{equation}
Thus, since we are in the space-periodic zero-average setting, the operator $\Lambda$ is well-defined as the multiplier by $2 \pi |k|$ in the Fourier space.
\\ \indent
The classical (fractional) Burgers equation corresponds to $f(u)=u^2/2$. The physical arguments justifying the small-scale estimates given above still hold for the class of functions $f$ considered here.
\\ \indent
There are two main types of physical motivations for studying the equation (\ref{fracBurgers}). In the field of nonlinear acoustics, it describes an asymptotic regime for detonation waves (see the paper of Clavin and Denet \cite{ClDe02} and also the introduction to the paper of Alfaro and Droniou \cite{AlDr12}). In the field of fluid mechanics, it has often been considered as a toy model for the SQG (Surface Quasi-Geostrophic) equation \cite{CTV15, CoCo04}.
\\ \indent
By the same heuristic arguments as above, denoting by $\ell$ the length scale for the solutions of (\ref{fracBurgers}), we obtain that $\ell^{-1} \sim \nu \ell^{-\alpha}$, and therefore $\ell \sim \nu^{\beta}$, where we denote by $\beta$ the quantity
\begin{equation} \label{beta}
\beta=\frac{1}{\alpha-1}.
\end{equation}
Note that for $\alpha=2$, i.e. for the classical Burgers equation, we have $\beta=1$. We see that the heuristic argument given above makes no sense for $\alpha \leq 1$, which suggests that the critical case is $\alpha=1$, where $\beta$ goes to $+\infty$. In the subcritical case $1<\alpha<2$, this dimensional analysis suggests that the results for our model should be the same as the results for the case $\alpha=2$ studied in \cite{BorD}, up to the replacement of $\nu$ by $\nu^{\beta}$. This is indeed the case: see Section~\ref{results}.
\\ \indent
Note that the arguments given above do not hold for $\alpha>2$: although the equation is still well-posed, since the operator $\Lambda^{\alpha}$ does not have a positive kernel, there is no good maximum principle for $u$.
\\ \indent
In the subcritical case $1<\alpha<2$, the well-posedness has first been proved for $x \in \R$ by Droniou, Gallou{\"e}t and Vovelle \cite{DGV03}; see also the earlier paper of Biler, Funaki and Woyczynski \cite{BiFuWo98} for partial results. Moreover in \cite{ADV07} Alibaud, Droniou and Vovelle proved that for large smooth initial data, the solutions are not necessarily smooth in the supercritical case $0<\alpha<1$.
\\ \indent
The well-posedness in the critical case $\alpha=1$ has been proved by Kiselev, Nazarov and Shterenberg \citep{KNS08} in the space-periodic setting using a modulus of continuity. This paper also contains a sketch of the proof of the well-posedness in the subcritical setting. For the sake of completeness, we include a detailed proof using the mild solution technique \cite{DZ92} in the Appendix. Note also that in \cite{KNS08} the supercritical ill-posedness result of \cite{ADV07} is extended to the space-periodic setting.
\\ \indent
An alternative proof of the well-posedness in the critical case has been given by Constantin and Vicol \cite{CoVi12} using a nonlinear maximum principle of the C{\'o}rdoba-C{\'o}rdoba type \cite{CoCo04}. A multidimensional generalisation of this result has been proved by Chan and Czubak \cite{ChCz10} extending the techniques of Caffarelli and Vasseur \cite{CaVa10}.
\\ \indent
The fractional Burgers equation has also been considered in a variety of other settings.
For a probabilistic interpretation of (\ref{fracBurgers}), see the papers of Jourdain, M{\'e}l{\'e}ard and Woyczynski \cite{JMW05}, Jourdain and Roux \cite{JoRo11} and Truman and Wu \cite{TrWu06}. For a proof of ergodicity for the fractional Burgers equation with space-time white noise, see the paper of Brzezniak, Debbi and Goldys \cite{BDG11}.

\subsection{Additional comments}
\smallskip
\indent
Estimating small-scale quantities for
\\
nonlinear PDEs is motivated by the problem of turbulence: for more information, see the book of Frisch \cite{Fri95} and the pioneering mathematically rigorous papers of Kuksin \cite{Kuk97GAFA,Kuk99GAFA}.
\\ \indent
In the same way as in \cite{BorD}, our estimates hold in average on a time interval $[T_1,T_2]$. In other words, we consider a time range during which we have the transitory behaviour which is referred to as decaying Burgers turbulence \cite{BK07}. This time interval depends only on $f$ and, through the quantity $D$ (see (\ref{D})), on $u_0$: thus it does not depend on $\nu$.
\\ \indent
When studying the typical behaviour for solutions of nonrandom PDEs, it is common to avoid pathological initial data, sometimes considering some type of averaging: see for instance \cite{BuTz08}. This is due to the lack of a random mechanism which allows to get solutions out of \enquote{bad} regions of the phase space, in particular in Hamiltonian PDEs. Here, the situation is more transparent: a non-zero initial condition $u_0$ is as generic as the ratio between the orders of $(u_0)_x$ and of $u_0$ itself. This ratio can be bounded from above using the quantity $D$:
\begin{equation} \label{D}
D=\max (|u_0|_{1}^{-1},\ | u_0 |_{1,\infty} )>1
\end{equation}
(see Section~\ref{Sob} for the notation).
\\ \indent
Note that for $m \in \lbrace 0,\ 1\rbrace$ and $1 \leq p \leq \infty$, we have:
\begin{equation} \label{Dp}
D^{-1} \leq |u_0|_{m,p} \leq D.
\end{equation}
The physical meaning of $D$ is that it gives a lower bound for the ratio between the amount of energy $\frac{1}{2} \int_{S^1}{u^2}$ initially contained in the system and its rate of dissipation $-\nu \int_{S^1}{u \Lambda^{\alpha} u}=\nu \Vert u \Vert_{\alpha/2}$ (see Section~\ref{Sob} for the notation).
\medskip
\\ \indent
In a future work we will consider the equation (\ref{fracBurgers}) with an additive random force, in a setting similar to \cite{BorK,BorW}. We expect to obtain the same results as in those papers up to the replacement of $\nu$ by $\nu^{\beta}$, i.e. with the same modifications as the ones in our paper with respect to \cite{BorD}.

\subsection{Plan of the paper} \label{plan}
\smallskip
 \indent
We introduce the notation and the setup in Section~\ref{nota}. We present the main results of our paper in Section~\ref{results}.
\\ \indent
In Section~\ref{upper}, we prove upper estimates for the positive and the negative parts of the quantity $\partial u/\partial x$. We use in a crucial way the nonlinear maximum principle of Constantin and Vicol \cite{CoVi12}. This result allows us to obtain upper bounds for the Sobolev norms $|u|_{m,p}$. In Section~\ref{lower}, using the results of Sections~\ref{upper}-\ref{turbupper}, we obtain time-averaged lower bounds for the Sobolev norms $|u|_{m,p}$. The upper and lower bounds are sharp, i.e. they coincide up to a $\nu$-independent multiplicative constant.
\\ \indent
In Sections~\ref{turbupper} and \ref{turblower} we obtain $\nu$-uniform sharp upper and lower bounds for the small-scale quantities corresponding to the flow $u(t,x)$, and we analyse the meaning of these results in terms of the theory of turbulence. Moreover, in Section~\ref{turbupper} we prove a crucial upper estimate for some fractional Sobolev norms, which will be used in Section~\ref{lower}.

\section{Notation and setup} \label{nota}
\smallskip
 \indent
\textbf{Agreement}: In the whole paper, all functions that we consider are real-valued and the space variable $x$ belongs to $S^1=\R/\Z$.

\subsection{Sobolev spaces} \label{Sob}
\smallskip
 \indent
Consider a zero mean value integrable function $v$ on $S^1$. 
For $p \in [1,\infty)$, we denote its $L_p$ norm
$$
\Bigg( \int_{S^1}{|v|^p} \Bigg)^{1/p}
$$
by $\left|v\right|_p$. The $L_{\infty}$ norm is by definition
$$
\left|v\right|_{\infty}=\esssup_{x \in S^1} |v(x)|.
$$
The $L_2$ norm is denoted by  $\left|v\right|$, and $\left\langle \cdot,\cdot\right\rangle$ stands for the $L_2$ scalar product. From now on $L_p,\ p \in [1,\infty],$ denotes the space of zero mean value functions in $L_p(S^1)$. Similarly, $C^{\infty}$ is the space of $C^{\infty}$-smooth zero mean value functions on $S^1$.
\\ \indent
For a nonnegative integer $m$ and $p \in [1,\infty]$, $W^{m,p}$ stands for the Sobolev space of zero mean value functions $v$ on $S^1$ with finite norm
\begin{equation} \nonumber
\left|v\right|_{m,p}=\left|\frac{d^m v}{dx^m}\right|_p.
\end{equation}
In particular, $W^{0,p}=L_p$ for $p \in [1,\infty]$. For $p=2$, we denote $W^{m,2}$ by $H^m$, and abbreviate the corresponding norm as $\left\|v\right\|_m$. 
\\ \indent
Note that since the length of $S^1$ is $1$ and the mean value of $v$ vanishes, we have:
\begin{equation} \label{Sobelement}
|v|_1 \leq |v|_{\infty} \leq |v|_{1,1} \leq |v|_{1,\infty} \leq \dots \leq |v|_{m,1} \leq |v|_{m,\infty} \leq \dots
\end{equation} 
We recall a version of the classical Gagliardo-Nirenberg inequality: cf. \cite[Appendix]{DG95}.
\begin{lemm} \label{GN}
For a smooth zero mean value function $v$ on $S^1$,
$$
\left|v\right|_{\chi,r} \leq C \left|v\right|^{\theta}_{m,p} \left|v\right|^{1-\theta}_{q},
$$
where $m>\chi$, and $r$ is determined by
$$
\frac{1}{r}=\chi-\theta \Big( m-\frac{1}{p} \Big)+(1-\theta)\frac{1}{q},
$$
under the assumption $\theta=\chi/m$ if $p=1$ or $p=\infty$, and $\chi/m \leq \theta < 1$ otherwise. The constant $C$ depends on $m,p,q,\chi,\theta$.
\end{lemm}
\indent
Subindices $t$ and $x$, which can be repeated, denote partial differentiation with respect to the corresponding variables. We denote by $v^{(m)}$ the $m$-th derivative of $v$ in the variable $x$. The function $v(t,\cdot)$ is abbreviated as $v(t)$.
\\ \indent
For any $s \geq 0$, $H^{s}$ stands for the Sobolev space of zero mean value functions $v$ on $S^1$ with finite norm
\begin{equation} \label{Hsdef}
\left\|v\right\|_{s}=(2 \pi)^{s} \Big( \sum_{k \in \Z}{|k|^{2s} |\hat{v}^k|^2} \Big)^{1/2},
\end{equation}
where $\hat{v}^k$ are the complex Fourier coefficients of $v(x)$. For an integer $s=m$, this norm coincides with the previously defined $H^m$ norm. For $s \in (0,1)$, $\left\|v\right\|_{s}$ is equivalent to the norm
\begin{equation} \label{Sobolevfrac}
\left\|v\right\|^{'}_{s}=\Bigg( \int_{S^1} \Big(\int_0^1 {\frac{|v(x+\ell)-v(x)|^2}{\ell^{2s+1}} d \ell} \Big) dx \Bigg)^{1/2}
\end{equation}
(see \cite{Ada75, Tay96}).
\\ \indent
H{\"o}lder's inequality yields the following well-known interpolation inequality:
\begin{equation} \label{Hsinterpol}
\left\|v\right\|_{s_2} \leq \left\|v\right\|_{s_1}^{\theta} \left\|v\right\|_{s_3}^{1-\theta},\ s_1 \leq s_2 \leq s_3,
\end{equation}
where:
$$
\theta=\frac{s_3-s_2}{s_3-s_1}.
$$
\indent
Lemma~\ref{GN} and (\ref{Hsinterpol}) yield the following inequalities, which will be used in the proof of Lemma~\ref{lmubuinfty}.

\begin{rmq} \label{expcount}
Here and below, when we prove upper estimates, we take advantage of the fact that we are in the one-dimensional setting, which allows to easily predict the powers for different Sobolev norms: in the inequalities, the \enquote{size} of $|\cdot|_{m,p}$ is $m-1/p$ and the \enquote{size} of $\Vert \cdot \Vert_s$ is $s-1/2$. For instance, the inequality:
$$
\Vert v \Vert \lesssim |v|_{1}^{2/3} \Vert v \Vert_1^ {1/3},
$$
which is a particular case of Lemma~\ref{GN}, can be predicted using the fact that:
$$
0-\frac{1}{2} = \frac{2}{3} \Big(0-1 \Big) + \frac{1}{3} \Big( 1-\frac{1}{2} \Big).
$$
Following this principle allows us to guess the right exponents, in particular in the three following lemmas; however, this principle cannot become a systematic rule. Indeed, there are restictions for admissible exponents in Lemma~\ref{GN}, and moreover this lemma does not allow us to estimate $H^s$ norms for noninteger values of $s$. Thus, to prove inequalities, we have to use Lemma~\ref{GN} and (\ref{Hsinterpol}) each time.
\end{rmq}

\begin{lemm} \label{fracinterpolHm}
For a smooth function $v$, we have:
\begin{equation} \label{fracinterpolHmstat}
\left\|v\right\|_{m} \overset{m}{\lesssim} \left|v\right|_{1,1}^{1-\theta} \left\|v\right\|_{m+\gamma}^{\theta},\ m \geq 2,\ 0 \leq \gamma \leq 1,
\end{equation}
where:
$$
\theta=\frac{2m-1}{2m+2 \gamma-1}.
$$
\end{lemm}

\textbf{Proof.} Using first (\ref{Hsinterpol}) and then Lemma~\ref{GN} for the function $v'$, we get:
\begin{align} \nonumber
 \left\|v\right\|_{m} & \leq \left\|v\right\|_{1}^{\gamma/(m+\gamma-1)} \left\|v\right\|_{m+\gamma}^{(m-1)/(m+\gamma-1)}
\\ \nonumber
& \overset{m}{\lesssim} \left|v\right|_{1,1}^{\gamma(2m-2)/(2m-1)(m+\gamma-1)} \left\|v\right\|_{m}^{\gamma/(2m-1)(m+\gamma-1)}
\\ \nonumber
& \times \left\|v\right\|_{m+\gamma}^{(m-1)/(m+\gamma-1)}.
\end{align}
Dividing by $\left\|v\right\|_{m}^{\gamma/(2m-1)(m+\gamma-1)}$ on both sides of the inequality, since
$$
1-\frac{\gamma}{(2m-1)(m+\gamma-1)}=\frac{(m-1)(2m+2 \gamma-1)}{(2m-1)(m+\gamma-1)},
$$
we obtain that:
\begin{align} \nonumber
& \left\|v\right\|_{m}^{(m-1)(2m+2 \gamma-1)/(2m-1)(m+\gamma-1)}
\\ \nonumber 
& \overset{m}{\lesssim} \left|v\right|_{1,1}^{\gamma(2m-2)/(2m-1)(m+\gamma-1)} \left\|v\right\|_{m+\gamma}^{(m-1)/(m+\gamma-1)},
\end{align}
which yields that:
$$
 \left\|v\right\|_{m} \overset{m}{\lesssim} \left|v\right|_{1,1}^{1-\theta} \left\|v\right\|_{m+\gamma}^{\theta}.\ \square
$$

\begin{lemm} \label{fracinterpol1infty}
We have:
\begin{equation} \label{fracinterpol1inftystat}
\left|v\right|_{1,\infty} \overset{m}{\lesssim} \left|v\right|_{1,1}^{1-\theta'} \left\|v\right\|_{m+\gamma}^{\theta'},\ m \geq 2,\ 0 \leq \gamma \leq 1,
\end{equation}
where:
$$
\theta'=\frac{2}{2m+2 \gamma-1}.
$$
\end{lemm}

\textbf{Proof.} This results follows immediately from the previous lemma after observing that by Lemma~\ref{GN} we have:
$$
\left|v\right|_{1,\infty} \overset{m}{\lesssim} \left|v\right|_{1,1}^{(2m-3)/(2m-1)} \left\|v\right\|_{m}^{2/(2m-1)}.\ \square
$$

\begin{lemm} \label{fracinterpol12}
We have:
\begin{equation} \label{fracinterpol12stat}
\left\| v \right\|_1 \overset{m}{\lesssim} \left|v\right|_{1,1}^{1-\theta''} \left\|v\right\|_{m+\gamma}^{\theta''},\ m \geq 2,\ 0 \leq \gamma \leq 1,
\end{equation}
where:
$$
\theta''=\frac{1}{2m+2 \gamma-1}.
$$
\end{lemm}

\textbf{Proof.} This results follows immediately from the previous lemma after observing that by H{\"o}lder's inequality we have:
$$
\left\| v \right\|_1 {\lesssim} \left|v\right|_{1,1}^{1/2} \left|v\right|_{1,\infty}^{1/2}.\ \square
$$

\subsection{Notation} \label{notasubsect}
\smallskip
 \indent
In this paper, we study asymptotic properties of solutions to (\ref{fracBurgers}) for small values of $\nu$, i.e. we suppose that
$$
0 < \nu \ll 1.
$$
We assume that $f$ is infinitely differentiable and satisfies (\ref{strconvex}). Moreover, we assume that $f$ and its derivatives satisfy:
\begin{equation} \label{poly}
\forall m \geq 0,\ \exists h \geq 0,\ C_m>0:\ |f^{(m)}(x)| \leq C_m (1+|x|)^h,\quad x \in \R,
\end{equation}
where $h=h(m)$ is a function such that $1 \leq h(1) < 2$ (the lower bound for $h(1)$ follows from (\ref{strconvex})). The usual Burgers equation corresponds to $f(x)=x^2/2$.
\\ \indent
We recall that we restrict ourselves to the case in which the initial condition $u_0:=u(0)$ has zero space average. Integrating by parts in space, one deduces that $u(t)$ satisfies (\ref{zero}) for all $t$. Furthermore, we assume that $u_0 \in C^{\infty}$. We also assume that we are not in the case $u_0 \equiv 0$, corresponding to the trivial solution $u(t,x) \equiv 0$. This ensures that the quantity $D$ (see (\ref{D})) is well-defined.
\smallskip
\\ \indent
Subindices $t$ and $x$, which can be repeated, denote partial differentiation with respect to the corresponding variables. We denote by $v^{(m)}$ the $m$-th derivative of $v$ in the variable $x$. For shortness, the function $v(t,\cdot)$ is denoted by $v(t)$.
\bigskip
\\ \indent
\textbf{Agreements:}\ From now on, all constants denoted by $C$ with sub- or superindexes are positive. Unless otherwise stated, they depend only on $f$ and on $D$. By $C(a_1,\dots,a_k)$ we denote constants which also depend on parameters $a_1,\dots,a_k$. By $X \overset{a_1,\dots,a_k}{\lesssim} Y$ we mean that $X \leq C(a_1,\dots,a_k) Y$. The notation $X \overset{a_1,\dots,a_k}{\sim} Y$ stands for
$$
Y \overset{a_1,\dots,a_k}{\lesssim} X \overset{a_1,\dots,a_k}{\lesssim} Y.
$$
In particular, $X \lesssim Y$ and $X \sim Y$ mean that $X \leq C Y$ and
\\
$C^{-1} Y \leq X \leq C Y$, respectively.
\\ \indent
All constants are independent of the viscosity $\nu$. We denote by $u=u(t,x)$ a solution of (\ref{fracBurgers}) for an initial condition $u_0$. A relation where the admissible values of $t$ (respectively, $x$) are not specified is assumed to hold for all $t \geq 0$, or $t>0$ of the relation contains $t^{-1}$ (respectively, all $x \in S^1$).
\\ \indent
The brackets $\lbrace \cdot \rbrace$ stand for the averaging in time over an interval $[T_1,T_2]$, where $T_1,T_2$ only depend on $f$ and on $D$ (see (\ref{T1T2}) for their definition.)
\\ \indent
For $m \geq 0$, $p \in [1,\infty]$, $\gamma(m,p)$ is by definition the quantity
\\
$\max(0,m-1/p)$.
\\ \indent
We use the notation $g^{-}=\max(-g,0)$ and $g^{+}=\max(g,0)$.

\subsection{Notation in Sections~\ref{turbupper} and \ref{turblower}} \label{notaturb}
\smallskip
\indent
In Sections~\ref{turbupper} and \ref{turblower}, we study analogues of quantities which are important for hydrodynamical turbulence. We consider quantities in the physical space (structure functions) as well as in the Fourier space (energy spectrum). We assume that $\nu \leq \nu_0$. The value of $\nu_0>0$ will be chosen in (\ref{nu0eq}).
\\
\smallskip
\indent
We define the non-empty and non-intersecting intervals
\begin{equation} \nonumber
J_1=(0,\ C_1 \nu^{\beta}];\ J_2=(C_1 \nu^{\beta},\ C_2];\ J_3=(C_2,\ 1]
\end{equation}
corresponding to the \textit{dissipation range}, the \textit{inertial range} and the \textit{energy range} from the Kolmogorov 1941 theory of turbulence \cite{Fri95}. For the definition of $\beta$, see (\ref{beta}).
\\ \indent
The quantities $S_{p}(\ell)$ denote the averaged moments of the increments in space for the flow $u(t,x)$:
$$
S_{p}(\ell)= \Bigg\{ \int_{S^1}{|u(t,x+\ell)-u(t,x)|^p dx} \Bigg\},\ p \geq 0,\ 0< \ell \leq 1.
$$
The quantity $S_{p}(\ell)$ is the \textit{structure function} of $p$-th order. The flatness, which measures spatial intermittency \cite{Fri95}, is defined by:
\begin{equation} \label{flatness}
F(\ell)=S_4(\ell)/S_2^2(\ell).
\end{equation}
Finally, for $k \geq 1$, we define the (layer-averaged) energy spectrum by
\begin{equation} \label{spectrum}
E(k)=\Bigg\{ \frac{\sum_{|n| \in [M^{-1}k,Mk]}{|\hat{u}(n)|^2}}{\sum_{|n| \in [M^{-1}k,Mk]}{1}} \Bigg\},
\end{equation}
where $M \geq 1$ is a constant which will be specified later (see the proof of \cite[Theorem 6.11]{BorD}).
\\ \indent
For more comments on the small-scale features of the solution, see the paper \cite{BorD}.

\section{Main results} \label {results}
\indent
In our paper, in Sections~\ref{upper} and \ref{lower}, we prove sharp upper and lower bounds for moments of Sobolev norms of $u$, generalising the results in \cite{BorD}. These results for Sobolev norms of solutions are summed up in Theorem~\ref{avoir}. Namely, for $m \in \lbrace 0,1 \rbrace$ and $p \in [1,\infty]$ or for $m \geq 2$ and $p \in (1,\infty]$ we have:
\begin{equation} \label{avoirresults}
\Big( \lbrace {\left|u(t)\right|_{m,p}^{\kappa}} \rbrace \Big)^{1/\kappa} \overset{m,p,\kappa}{\sim} \nu^{-\beta \gamma},\quad \kappa>0,
\end{equation}
and on the other hand:
\begin{equation} \label{E}
\Big( \lbrace {\left\|u(t)\right\|_{s}^{\kappa}} \rbrace \Big)^{1/\kappa} \overset{s,\kappa}{\sim} \nu^{-\beta (s-1/2)},\quad s>1/2.
\end{equation}
We recall that by definition, $\gamma(m,p)=\max(0,m-1/p)$, and the brackets $\lbrace \cdot \rbrace$ stand for the averaging in time over an interval $[T_1,T_2]$ ($T_1, T_2$ only depend on $f$ and, through $D$, on $u_0$: see (\ref{T1T2})).
\\ \indent
In Section~\ref{turbupper} and \ref{turblower} we obtain sharp estimates for analogues of quantities characterising hydrodynamical turbulence. In what follows, we assume that $\nu \in (0,\nu_0]$, where $\nu_0 \in (0,1]$ depends only on $f$ and on $D$.
\\ \indent
First, as a consequence of (\ref{avoirresults})-(\ref{E}), we prove Theorem~\ref{avoir2}, which states that for $\ell \in J_1$:
$$
\quad \ \ \ S_{p}(\ell) \overset{p}{\sim} \left\lbrace \begin{aligned} & \ell^{p},\ 0 \leq p \leq 1. \\ & \ell^{p} \nu^{-\beta (p-1)},\ p \geq 1, \end{aligned} \right.$$
and for $\ell \in J_2$:
$$
S_{p}(\ell) \overset{p}{\sim} \left\lbrace \begin{aligned} & \ell^{p},\ 0 \leq p \leq 1. \\ & \ell,\ p \geq 1. \end{aligned} \right.
$$
Consequently, for $\ell \in J_2$ the flatness satisfies the estimate:
$$
F(\ell)=S_4(\ell)/S_2^2(\ell) \sim \ell^{-1}.
$$
This gives a rigorous proof of the fact that $u$ is highly intermittent in the inertial range. 
\\ \indent
Finally we get estimates for the spectral asymptotics of the decaying Burgulence. On one hand, for $m \geq 1$ we have:
$$
\lbrace |\hat{u}(k)|^2 \rbrace \overset{m}{\lesssim} k^{-2 m} {\Vert u \Vert_m^2 } \overset{m}{\lesssim} (k \nu^ {\beta} )^{-2m} \nu^{\beta}.
$$
In particular, $\lbrace |\hat{u}(k)|^2 \rbrace$ decreases at a faster-than-algebraic rate for $|k| \succeq \nu^{-\beta}$. On the other hand, by Theorem~\ref{spectrinert}, for $k$ such that $k^{-1} \in J_2$ the energy spectrum $E(k)$ satisfies:
$$
E(k) \sim k^{-2},
$$
where the quantity $M \geq 1$ in the definition of $E(k)$ depends only on $f$ and on $D$. 

\begin{rmq}
The main results of our paper are word-to-word the same as in the paper \cite{BorD} with $\nu$ replaced by $\nu^{\beta}$ and some modifications in the defintions of the constants (in particular the interval $[T_1,T_2]$).
\end{rmq}

\section{Upper estimates for Sobolev norms} \label{upper}
\smallskip
 \indent
We recall that $u=u(t,x)$ denotes a solution of (\ref{fracBurgers}) for an initial condition $u_0$. For more information on the notation, see Section~\ref{nota}.
\\ \indent
We begin by proving a key upper estimate for $u_x$. This estimate is well-known in the case $\alpha=2$ \cite{Kru64}.

\begin{lemm} \label{uxpos}
We have:
$$
u_x(t,x) \leq \min (D,\sigma^{-1} t^{-1} ).
$$
\end{lemm}

\textbf{Proof.}
Differentiating the equation (\ref{fracBurgers}) once in space, multiplying by $t$ and considering the function $v=tu_x$ we obtain that:
\begin{equation} \label{uxposeq}
v_t+t^{-1}(-v+f''(u) v^2)+f'(u)v_x=-\nu \Lambda^{\alpha} v.
\end{equation}
Now suppose that we are not in the trivial case $v \equiv 0$ and consider a point $(t_1,x_1)$ where $v$ reaches its maximum on the cylinder $S=[0,t] \times S^1$. We assume that $t_1>0$ and we remark that this maximum is positive since for every $t$, $v(t,\cdot)$ has zero space average. At $(t_1,x_1)$, Taylor's formula implies that we have $v_t \geq 0$ (since $t_1>0$) and $v_x=0$ (since $S^1$ has no boundary). Moreover, the positivity of the kernel associated to the operator $-\Lambda^{\alpha}$ for $\alpha \leq 2$ implies that $\Lambda^{\alpha} v \geq 0$ (cf. \cite{CoVi12}). Therefore, (\ref{uxposeq}) yields that:
$$
-v+f''(u) v^2 \leq 0,
$$
which implies that $v(t_1,x_1) \leq \sigma^{-1}$. Therefore we have $u_x(t,x) \leq \sigma^{-1} t^{-1}$ for all $t,x$.
\\ \indent
To prove that for all $x$, $u_x(t,x) \leq \min_{y \in S^1}{u_0(y)} \leq D$, we use a simpler version of the same argument applied to the function $u_x$. $\square$
\bigskip
\\ \indent
Since the space averages of $u(t)$ and $u_x(t)$ vanish, we get the following upper estimates. First, for $1 \leq p \leq +\infty$, we get:
\begin{align} \label{Lpupper}
&\left|u(t)\right|_{p} \leq \left|u(t)\right|_{\infty} \leq \int_{S^1}{u_x^{+}(t)} \leq \min (D, \sigma^{-1} t^{-1} ).
\end{align}
Then we get the following crucial estimate:
\begin{align} \nonumber
\left|u(t)\right|_{1,1}&=\int_{S^1}{u_x^{+}(t)}+\int_{S^1}{u_x^{-}(t)}=2 \int_{S^1}{u_x^{+}(t)}
\\ \label{W11}
& \leq 2\min (D, \sigma^{-1} t^{-1}).
\end{align}

\begin{lemm} \label{nonlin}
We have:
$$
\left|u(t)\right|_{1,\infty} {\lesssim}\ \nu^{- \beta}.
$$
\end{lemm}

\textbf{Proof.} 
By Lemma~\ref{uxpos}, it suffices to prove that we have:
$$
h=-u_x \lesssim \nu^{-\beta}.
$$
Differentiating the equation (\ref{fracBurgers}) once in space we obtain that:
\begin{equation} \label{uxeq}
h_t-f''(u) h^2+f'(u) h_x=-\nu \Lambda^{\alpha} h.
\end{equation}
Now consider a point $(t_1,x_1)$ where $h$ reaches its maximum on the cylinder $S=[0,t] \times S^1$ and suppose that we are not in the trivial case
$h \equiv 0$. Suppose also that $t_1>0$ (else $\max_{(t,x) \in S}{|h(t,x)|} \leq D$) and denote this maximum by $M$.
Then at $(t_1,x_1)$ by Taylor's formula we have $h_t \geq 0$ and $h_x=0$. On the other hand, by \cite[Theorem 2.3]{CoVi12} we have one of the two following situations:
\\ \indent
$a)$:\ $h(t_1,x_1) \lesssim \max_{}(-u)$, and therefore by (\ref{Lpupper}), $h(t_1,x_1) \lesssim 1$.
\\ \indent
$b)$:\ $\Lambda^{\alpha} h(t_1,x_1) \gtrsim h^{1+\alpha}(t_1,x_1)/|u(t_1,\cdot)|^{\alpha}_{\infty}$, and therefore by (\ref{Lpupper}),\\ $\Lambda^{\alpha} h(t_1,x_1) \gtrsim M^{1+\alpha}$.
\\ \indent
In the situation $b)$, (\ref{uxeq}) yields that at the point $(t_1,x_1)$:
$$
-f''(u(t_1,x_1)) M^2 \lesssim -\nu M^{1+\alpha},
$$
and therefore by (\ref{strconvex}) we get:
$$
M \lesssim \nu^{-1/(\alpha-1)}=\nu^{-\beta}.\  \square
$$

\begin{lemm} \label{upper1}
We have the inequality
$$
\left\|u(t)\right\|^2_{1} \lesssim \nu^{-\beta}.
$$
\end{lemm}

\textbf{Proof.} It suffices to use (\ref{W11}) and Lemma~\ref{nonlin} and to apply H{\"o}lder's inequality.  $\square$
\medskip
\\ \indent
Now we prove an important auxiliary lemma, which plays the same role as \cite[Lemma 5.2.]{BorD}. The modifications in the exponents which follow from the modifications of the Sobolev norms can be guessed using Remark~\ref{expcount}.

\begin{lemm} \label{lmubuinfty}
For every $m \geq 1$ there exist $C_m>0$ and a natural number $n'=n'(m)$ such that for $v \in C^{\infty}$,
\begin{align} \label{polyestimate}
N_m(v) &:=\left| \left\langle v^{(m)}, (f(v))^{(m+1)} \right\rangle\right|
\\ \nonumber
& \leq C_m (1+\left|v\right|_{1,1} )^{n'} \left\|v\right\|_{m+\alpha/2}^{4m/(2m+\alpha-1)}.
\end{align}
\end{lemm}

\textbf{Proof.}
Fix $m \geq 1$. We denote $\left|v\right|_{1,1}$ by $N$ and we recall that $\left|v\right|_{\infty} \leq N$. Let $C'$ denote various expressions of the form $C_m (1+N)^{n(m)}$. Integrating by parts, we get:
\begin{align} \nonumber
& N_m(v) =  \left|\left\langle v^{(2m)}, (f(v))^{(1)} \right\rangle\right|= \left|\left\langle v^{(m)}, (f(v))^{(m+1)} \right\rangle\right|
\\ \nonumber
& \leq  C(m) \sum_{k=1}^m\ \sum_{\substack{1 \leq a_1 \leq \dots \leq a_k \leq m-1 \\ a_1+ \dots+a_k = m+1}} \int_{S^1}{\left| v^{(m)} v^{(a_1)} \dots v^{(a_k)} f^{(k)}(v) \right|}
\\ \nonumber
& + C(m) \left| \int_{S^1}{ (v^{(m)})^ 2 v' f''(v) } \right| + \left| \int_{S^1}{ v^{(m)} v^{(m+1)} f'(v) } \right|
\\ \nonumber
& =  C(m) \sum_{k=1}^m\ \sum_{\substack{1 \leq a_1 \leq \dots \leq a_k \leq m-1 \\ a_1+ \dots+a_k = m+1}} \int_{S^1}{\left| v^{(m)} v^{(a_1)} \dots v^{(a_k)} f^{(k)}(v) \right|}
\\ \nonumber
& + C(m) \left|\frac{3}{2} \int_{S^1}{ f''(v) (v^{(m)})^2 v' } \right|
\\ \nonumber
& \leq  C(m) \max_{x \in [-N,N]}\ \max(f'(x),\dots f^{(m)}(x))
\\ \nonumber
& \times \sum_{k=1}^m\ \sum_{\substack{1 \leq a_1 \leq \dots \leq a_k \leq m-1 \\ a_1+ \dots+a_k = m+1}} \int_{S^1} | v^{(a_1)} \dots v^{(a_k)} v^{(m)} |
\\ \nonumber
& + C(m) \left| \int_{S^1}{ f''(v) (v^{(m)})^2 v' } \right|.
\end{align}
Therefore by (\ref{poly}) and (\ref{Lpupper}) we get:
\begin{align} \nonumber
N_m(v) \leq & \ \underbrace{C' \sum_{\substack{1 \leq a_1 \leq \dots \leq a_k \leq m-1 \\ a_1+ \dots+a_k = m+1}} \int_{S^1} { | v^{(a_1)} \dots v^{(a_k)} v^{(m)} |}}_{N_{m,1}}
\\ \nonumber
& +\underbrace{C' \int_{S^1}{ \left|(v^{(m)})^2 v' \right| } }_{N_{m,2}}.
\end{align}
By (\ref{upper}), Lemma~\ref{fracinterpolHm} and Lemma~\ref{fracinterpol1infty}, there exists a constant $\epsilon(m)>0$ such that:
\begin{align} \nonumber
N_{m,2} & \leq C' |v|_{1,\infty} \left\|v\right\|_{m}^{2}
\\ \nonumber
& \leq C' |v|_{1,1}^{\epsilon(m)} \left\|v\right\|_{m+\alpha/2}^{2/(2m+\alpha-1)} \left\|v\right\|_{m+\alpha/2}^{(4m-2)/(2m+\alpha-1)}
\\ \nonumber
& \leq C' \left\|v\right\|_{m+\alpha/2}^{4m/(2m+\alpha-1)}.
\end{align}
Now it remains to deal with $N_{m,1}$. Using first H{\"o}lder's inequality, then Lemma~\ref{GN} and finally (\ref{Hsinterpol}), we get:
\begin{align} \nonumber
& N_{m,1} \leq C' \sum_{\substack{1 \leq a_1 \leq \dots \leq a_k \leq m-1 \\ a_1+ \dots+a_k = m+1}}{| v^{(a_1)}|_{\infty} \dots |v^{(a_{k-1})}|_{\infty} \left\|v\right\|_{a_k} \left\|v\right\|_{m}}
\\ \nonumber
& \leq C' \left\|v\right\|_{m} \sum_{\substack{1 \leq a_1 \leq \dots \leq a_k \leq m-1 \\ a_1+ \dots+a_k = m+1}}{ (\left\| v \right\|_{a_1}^{1/2} \left\| v \right\|_{a_1+1}^{1/2}) \dots (\left\|v\right\|_{a_{k-1}}^{1/2} \left\| v\right\|_{a_{k-1}+1}^{1/2})  \left\| v\right\|_{a_k}}
\\ \nonumber
& \leq C' \left\|v\right\|_{m} \sum_{\substack{1 \leq a_1 \leq \dots \leq a_k \leq m-1 \\ a_1+ \dots+a_k = m+1}}{ \Big[ (\left\| v \right\|_{1}^{(2m-2a_1-1)/2(m-1)} \left\| v \right\|_{m}^{(2a_1-1)/2(m-1)}) \dots}
\\ \nonumber
& \times (\left\|v\right\|_{1}^{(2m-2a_{k-1}-1)/2(m-1)} \left\| v\right\|_{m}^{(2a_{k-1}-1)/2(m-1)})
\\ \nonumber
& \times (\left\|v\right\|_{1}^{(2m-2a_{k})/2(m-1)} \left\| v\right\|_{m}^{(2a_{k}-2)/2(m-1)}) \Big]
\\ \nonumber
& \leq C' \left\|v\right\|_{m} \sum_{k=1}^{m+1}{ \left\| v \right\|_{1}^{(2mk-2m-k-1)/2(m-1)} \left\| v \right\|_{m}^{(2m-k+1)/2(m-1)}}.
\end{align}
Using Lemma~\ref{fracinterpolHm} and Lemma~\ref{fracinterpol12}, we get that for every $k$ there exists $\epsilon'(k)$ such that:
\begin{align} \nonumber
& N_{m,1} \leq C' \left\|v\right\|_{m} \sum_{k=1}^{m+1}{C' \left| v \right|_{1,1}^{\epsilon'(k)} \left\| v \right\|_{m+\alpha/2}^{(2m+1)/(2m+\alpha-1)}}.
\end{align}
Using once again Lemma~\ref{fracinterpolHm}, we get:
\begin{align} \nonumber
& N_{m,1} \leq C' \left\|v\right\|_{m} \left\| v \right\|_{m+\alpha/2}^{(2m+1)/(2m+\alpha-1)}
\\ \nonumber
& \leq C' \left\| v \right\|_{m+\alpha/2}^{4m/(2m+\alpha-1)}.\ \square
\end{align}
\medskip \\ \indent
The following result shows that there is a strong nonlinear damping which prevents the successive derivatives of $u$ from becoming too large.

\begin{lemm} \label{uppermaux}
For integer values of $m \geq 1$,
$$
\left\|u(t)\right\|^{2}_{m} \overset{m}{\lesssim} \max (\nu^{-(2m-1) \beta}, t^{-(2m-1)}).
$$
(we recall that $\beta=1/(\alpha-1)$).
\end{lemm}

\textbf{Proof.} Fix $m \geq 1$. Denote
$$
x(t)=\left\|u(t)\right\|^{2}_{m}.
$$
We claim that the following implication holds:
\begin{align} \label{decrm}
&x(t) \geq \bar{C} \nu^{-(2m-1) \beta} \Longrightarrow \frac{d}{dt} x(t) \leq -(2m-1) x(t)^{2m/(2m-1)},
\end{align}
where $\bar{C}$ is a fixed positive number, chosen later. Below, all constants denoted by $C$ do not depend on $\bar{C}$.
\\ \indent
Indeed, assume that $ x(t) \geq \bar{C} \nu^{-(2m-1) \beta}.$ Now denote
$$
y(t)=\left\|u(t)\right\|^{2}_{m+\alpha/2}.
$$
By Lemma~\ref{fracinterpolHm} and (\ref{W11}) we get:
\begin{align} \label{auxpower}
y(t) & \gtrsim x(t)^{(2m+\alpha-1)/(2m-1)} 
\\ \label{auxpower2}
& \geq \bar{C}^{(2m+\alpha-1)/(2m-1)} \nu^{-(2m+\alpha-1) \beta}.
\end{align}
On the other hand, integrating by parts in space and using Lemma~\ref{lmubuinfty}, we get the following energy dissipation relation:
\begin{align} \nonumber
\frac{d}{dt}  x(t) & = - 2 \nu  \left\|u(t)\right\|_{m+\alpha/2}^2-2\left\langle u^{(m)}(t), (f(u(t)))^{(m+1)}\right\rangle
\\ \indent
& \leq - 2 \nu  \left\|u(t)\right\|_{m+\alpha/2}^2 + C \left\| u(t) \right\|_{m+\alpha/2}^{4m/(2m+\alpha-1)}.
\\ \indent
& = (- 2 \nu y(t)^{1/\beta(2m+\alpha-1)}+C) y(t)^{2m/(2m+\alpha-1)}.
\end{align}
Thus, using (\ref{auxpower}) and (\ref{auxpower2}), we get that for $\bar{C}$ large enough:
\begin{align} \nonumber
\frac{d}{dt}  x(t) &\leq (-C \bar{C}^{1/\beta(2m-1)} + C ) x(t)^{2m/(2m-1)}.
\end{align}
Thus we can choose $\bar{C}$ in such a way that the implication (\ref{decrm}) holds. We claim that
\begin{equation} \label{decrmcor}
x(t) \leq \max ( \bar{C} \nu^{-(2m-1) \beta}, t^{-(2m-1)} ).
\end{equation}
Indeed, if $x(s) \leq \bar{C} \nu^{-(2m-1) \beta}$ for some $s \in \left[0,t\right]$, then the assertion (\ref{decrm}) ensures that $x(s)$ remains below this threshold up to time $t$.
\\ \indent
Now, assume that $x(s) > \bar{C} \nu^{-(2m-1) \beta}$ for all $s \in \left[0,t\right]$. Denote
$$
\tilde{x}(s)=(x(s))^{-1/(2m-1)},\ s \in \left[0,t\right].
$$
By (\ref{decrm}) we get $d\tilde{x}(s)/ds \geq 1$. Therefore $\tilde{x}(t) \geq t$ and $x(t) \leq t^{-(2m-1)}$. Thus in this case, the inequality (\ref{decrmcor}) still holds. This proves the lemma's assertion. $\square$
\bigskip \\ \indent
Applying the inequality (\ref{Hsinterpol}) we get the following result:
\begin{lemm} \label{uppermauxcor}
For $s \geq 1$, $s$ not necessarily being an integer,
$$
\left\|u(t)\right\|^{2}_{s} \overset{m}{\lesssim} \max (\nu^{-(2s-1) \beta}, t^{-(2s-1)}).
$$
\end{lemm}

The proof of the following lemma is word-to-word the same as the proof of \cite[Lemma 5.4]{BorD}.

\begin{lemm} \label{upperwmp}
For $m \in \lbrace 0,1 \rbrace$ and $p \in [1,\infty]$, or for $m \geq 2$ and $p \in (1,\infty]$ we have:
\begin{align} \nonumber
\left|u(t)\right|_{m,p} & \overset{m,p}{\lesssim} \max(\nu^{-\gamma \beta},t^{-\gamma}),
\end{align}
where we denote
$$
\gamma=\max(0,\ m-1/p).
$$
\end{lemm}

\section{Upper estimates for small-scale quantities} \label{turbupper}
\smallskip
\indent
In this section, we study analogues of quantities which are important for the study of hydrodynamical turbulence. For notation for these quantities and the ranges $J_1,\ J_2,\ J_3$, see Section~\ref{notaturb}. The statements and the proofs are word-to-word the same as in the case $\alpha=2$, up to the replacement of $\nu$ by $\nu^{\beta}$. Therefore we will omit the proofs.
\\ \indent
Moreover, in this section, we prove an upper estimate for the norms $\Vert u \Vert_{s}$,\ $s \in (1/2,\ 1)$, which will play a crucial role for the lower estimates in Section~\ref{lower}.

\begin{lemm} \label{upperdiss}
For $\ell \in [0,1]$,
$$
S_{p}(\ell) \overset{p}{\lesssim} \left\lbrace \begin{aligned} & \ell^{p},\ 0 \leq p \leq 1. \\ & \ell^{p} \nu^{-\beta (p-1)},\ p \geq 1. \end{aligned} \right.
$$
\end{lemm}

\begin{lemm} \label{upperinert}
For $\ell \in J_2 \cup J_3$,
$$
S_{p}(\ell) \overset{p}{\lesssim} \left\lbrace \begin{aligned} & \ell^{p},\ 0 \leq p \leq 1. \\ & \ell,\ p \geq 1. \end{aligned} \right.
$$
\end{lemm}

\begin{lemm} \label{Hsupper}
We have
$$
\Vert u(t) \Vert_{s}^2 \lesssim \nu^{-\beta(2s-1)}, s \in (1/2,\ 1).
$$
\end{lemm}

\textbf{Proof.} This proof follows the lines of \cite[Lemma 4.12]{BorW}.
\\ \indent
By (\ref{Sobolevfrac}) we have:
\begin{align} \nonumber
\Vert u(t) \Vert_{s}^2 \lesssim \int_{S^1} \Big(\int_0^1 {\frac{|u(t,x+\ell)-u(t,x)|^2}{\ell^{(2s+1)}} d \ell} \Big) dx.
\end{align}
Consequently, by Fubini's theorem,
\begin{align} \nonumber
\Vert u(t) \Vert_{s}^2 \rbrace &\lesssim \int_0^1 \frac{1}{\ell^{(2s+1)}} \Big( \int_{S^1}{|u(t,x+\ell)-u(t,x)|^2 dx} \Big) d \ell 
\\ \nonumber
&= \int_{0}^{1}{\frac{S_2(\ell)}{\ell^{(2s+1)}} d\ell}
=\int_{J_1}{\frac{S_2(\ell)}{\ell^{(2s+1)}} d\ell}+\int_{J_2}{\frac{S_2(\ell)}{\ell^{(2s+1)}} d\ell}+\int_{J_3}{\frac{S_2(\ell)}{\ell^{(2s+1)}} d\ell}.
\end{align}
By Lemma~\ref{upperdiss} we get:
$$
\int_{J_1}{\frac{S_2(\ell)}{\ell^{(2s+1)}} d\ell} \lesssim \int_{0}^{C_1 \nu^{\beta}}{\frac{\ell^2 \nu^{-\beta}}{\ell^{(2s+1)}} d\ell} \sim \nu^{-\beta} \nu^{\beta(2-2s)} = \nu^{-\beta(2s-1)} 
$$
and
$$
\int_{J_2}{\frac{S_2(\ell)}{\ell^{(2s+1)}} d\ell} \lesssim \int_{C_1 \nu^{\beta}}^{C_2}{\frac{\ell}{\ell^{(2s+1)}} d\ell} \sim \nu^{-\beta(2s-1)}.
$$
Finally, by (\ref{Lpupper}) we get:
$$
\int_{J_3}{\frac{S_2(\ell)}{\ell^{(2s+1)}} d\ell} \leq C C_2^{-(2s+1)} \leq C.
$$
Thus,
$$
\Vert u\Vert_{s}^2 \lesssim \nu^{-\beta(2s-1)}.\ \qed
$$

\section{Lower estimates for Sobolev norms} \label{lower}
\smallskip
 \indent
We define
\begin{equation} \label{T1T2}
T_1=\frac{1}{4}D^{-2} \tilde{C}^{-1};\quad T_2=\max \Bigg(\ \frac{3}{2} T_1,\quad 2D\sigma^{-1} \Bigg),
\end{equation}
where $\tilde{C}$ is a constant such that for all $t$, $\left\|u(t)\right\|_{\alpha/2}^2 \leq \tilde{C} \nu^{-1}$ (cf.
\\
Lemma~\ref{Hsupper}). Note that $T_1$ and $T_2$ do not depend on the viscosity coefficient $\nu$. 
\\ \indent
From now on, for any function $A(t)$, $\lbrace A(t) \rbrace$ is by definition the time average
$$
\lbrace A(t) \rbrace=\frac{1}{T_2-T_1} \int_{T_1}^{T_2}{A(s)\ ds}.
$$
\\ \indent
The first quantities that we estimate from below are the Sobolev norms $\lbrace |u(t)|_p^2 \rbrace,\ p \in [1,\infty]$.

\begin{lemm} \label{finitetimep}
For $p \in [1,\infty]$, we have:
$$
\lbrace |u(t)|_p^2 \rbrace \gtrsim 1.
$$
\end{lemm}

\textbf{Proof.} It suffices to prove the lemma's statement for $p=1$. But this case follows from the case $p=2$. Indeed, by H{\"o}lder's inequality and (\ref{Lpupper}) we get:
$$
\lbrace |u(t)|_1^2 \rbrace \geq \lbrace |u(t)|_{\infty}^{-2} |u(t)|^{4} \rbrace \gtrsim \lbrace |u(t)|^{4} \rbrace \geq \lbrace |u(t)|^{2} \rbrace^2.
$$
\indent
Integrating by parts in space, we get the dissipation identity
\begin{align} \label{dissip}
\frac{d}{dt} \left| u(t)\right|^2 &= \int_{S^1}{(-2 uf'(u)u_x+2 \nu u \Lambda^{\alpha} u)}=-2 \nu \left\| u(t)\right\|_{\alpha/2}^{2}.
\end{align}
Thus, integrating in time and using (\ref{D}) and Lemma~\ref{Hsupper}, we obtain that for $t \in [0,\ 3T_1/2]$ we have the following uniform lower bound:
\begin{align} \label{dissipaux}
|u(t)|^2 &= |u_0|^2 - 2 \nu \int_{0}^{t}{\left\|u(t)\right\|_{\alpha/2}^2} \geq D^{-2} - 3 T_1 \tilde{C} = D^{-2}/4.
\end{align}
Thus,
$$
\lbrace |u(t)|^2 \rbrace \geq  \frac{1}{T_2-T_1} \int_{T_1}^{3T_1/2}{|u(t)|^2} \geq  \frac{D^{-2} T_1}{8(T_2-T_1)}.\ \square
$$
\medskip
\\ \indent
Now we prove a key estimate for $\lbrace \left\|u(t)\right\|_{\alpha/2}^2 \rbrace$.

\begin{lemm} \label{finitetime}
We have
$$
\lbrace \left\|u(t)\right\|_{\alpha/2}^2 \rbrace \gtrsim \nu^{-1}.
$$
\end{lemm}

\textbf{Proof.} In the same way as in (\ref{dissipaux}), we prove that $|u(T_1)|^2 \geq D^{-2}/2$. Thus, using (\ref{Lpupper}) ($p=2$) and integrating in time the equality (\ref{dissip}) we get:
\begin{align} \nonumber
\lbrace \left\|u(t)\right\|_{\alpha/2}^2 \rbrace &= \frac{1}{2 \nu (T_2-T_1)} (|u(T_1)|^2-|u(T_2)|^2)
\\ \nonumber
& \geq \frac{1}{2 \nu (T_2-T_1)} \Big(\frac{1}{2}D^{-2}-\sigma^{-2} T_2^{-2} \Big)
\\ \nonumber
& \geq \frac{D^{-2}}{8 (T_2-T_1) } \nu^{-1},
\end{align}
which proves the lemma's assertion.\ $\square$

\begin{cor} \label{finitetimeH1}
We have
$$
\lbrace \left\|u(t)\right\|_{1}^2 \rbrace \gtrsim \nu^{-\beta}.
$$
\end{cor}

\textbf{Proof.} By (\ref{Hsinterpol}) and H{\"o}lder's inequality we get:
$$
\lbrace \left\|u(t)\right\|_{\alpha/2}^{2} \rbrace \leq \lbrace \left\|u(t)\right\|_{(1+\alpha)/4}^{2} \rbrace^{(4-2 \alpha)/(3-\alpha)} \lbrace \left\|u(t)\right\|_{1}^2 \rbrace^{(\alpha-1)/(3-\alpha)},
$$
and therefore
$$
\lbrace \left\|u(t)\right\|_{1}^2 \rbrace \geq \lbrace \left\|u(t)\right\|_{\alpha/2}^{2} \rbrace^{(3-\alpha)/(\alpha-1)} \lbrace \left\|u(t)\right\|_{(1+\alpha)/4}^{2} \rbrace^{-(4-2\alpha)/(\alpha-1)}.
$$
Thus, by Lemma~\ref{Hsupper} and Lemma~\ref{finitetime} we get:
$$
\lbrace \left\|u(t)\right\|_{1}^2  \gtrsim \nu^{-(3-\alpha)/(\alpha-1)} \nu^{(2-\alpha)/(\alpha-1)}=\nu^{-\beta}.\ \square
$$
\smallskip \\ \indent
This time-averaged lower bound yields similar bounds for other
\\
Sobolev norms: the proofs are word-to-word the same as in \cite{BorD}. We recall that $\gamma=m-1/p$. The result in Lemma~\ref{finalexp} can easily be extended to fractional Sobolev norms.

\begin{lemm} \label{finalexp}
For $m \geq 2$,
$$
\lbrace \left\|u(t)\right\|_m^2 \rbrace \overset{m}{\gtrsim} \nu^{-(2m-1) \beta}.
$$
\end{lemm} 

\begin{lemm} \label{finalexpbis}
For $m \geq 0$ and $p \in [1,\infty]$,
$$
\lbrace \left|u(t)\right|_{m,p}^2 \rbrace^{1/2} \overset{m,p}{\gtrsim} \nu^{-\gamma}.
$$
\end{lemm}

\begin{lemm} \label{finalexpter}
For $m \geq 0$ and $p \in [1,\infty]$,
$$
\lbrace \left|u(t)\right|_{m,p}^{\kappa} \rbrace^{1/\kappa} \overset{m,p,\kappa}{\gtrsim} \nu^{-\gamma},\quad \kappa>0.
$$
\end{lemm}

The following theorem sums up the main results of Sections~\ref{upper} and \ref{lower}.

\begin{theo} \label{avoir}
For $m \in \lbrace 0,1 \rbrace$ and $p \in [1,\infty]$, or for $m \geq 2$ and $p \in (1,\infty]$ we have:
\begin{equation} \label{avoirsim}
\Big( \lbrace \left|u(t)\right|_{m,p}^{\kappa} \rbrace \Big)^{1/\kappa} \overset{m,p,\kappa}{\sim} \nu^{-\beta \gamma},\qquad \kappa>0,
\end{equation}
and for $s>1/2$ we have:
\begin{equation} \label{avoirsim2}
\Big( \lbrace \left\|u(t)\right\|_{s}^{\kappa} \rbrace \Big)^{1/\kappa} \overset{s,\kappa}{\sim} \nu^{-\beta (s-1/2)},\qquad \kappa>1/2,
\end{equation}
where $\lbrace \cdot \rbrace$ denotes time-averaging over $[T_1,T_2]$. The upper estimates in (\ref{avoirsim}) hold without time-averaging, uniformly for $t$ separated from $0$. Namely, we have:
$$
\left|u(t)\right|_{m,p} \overset{m,p}{\lesssim} \max(t^{-\beta \gamma},\ \nu^{-\gamma}).
$$
On the other hand, the lower estimates hold for all $m \geq 0$ and $p \in [1,\infty]$, and also for $s=1/2$. 
\end{theo}

\textbf{Proof.} 
Upper estimates follow from Lemma~\ref{upperwmp}, and lower estimates from Lemma~\ref{finalexpter}.\ $\square$

\section{Lower estimates for small-scale quantities} \label{turblower}
\smallskip
\indent
In this section, we study analogues of quantities which are important for the study of hydrodynamical turbulence. For notation for these quantities and the ranges $J_1,\ J_2,\ J_3$, see Section~\ref{notaturb}. The statements and the proofs are word-to-word the same as in the case $\alpha=2$, up to the replacement of $\nu$ by $\nu^{\beta}$. Therefore we will omit the proofs.
\\ \indent
Provided $\nu \leq \nu_0$, all estimates hold independently of the viscosity $\nu$. We recall that the brackets $\lbrace \cdot \rbrace$ stand for the averaging in time over an interval $[T_1,T_2]$: see (\ref{T1T2}).

\begin{defi}
For $K>1$, we denote by $L_K$ the set of all $t \in [T_1,T_2]$ such that the assumptions
\begin{align} \label{condi}
& K^{-1} \leq |u(t)|_{\infty} \leq  \max u_x(t) \leq K
\\ \label{condii}
& K^{-1} \nu^{-\beta} \leq  |u(t)|_{1,\infty} \leq K \nu^{-\beta}
\\ \label{condiii}
& |u(t)|_{2,\infty} \leq K \nu^{-2 \beta}
\end{align}
hold.
\end{defi}

\begin{lemm} \label{typical}
There exist constants $C,K_1>0$ such that for $K \geq K_1$, the Lebesgue measure of $L_K$ satisfies $\lambda(L_K) \geq C$.
\end{lemm}

Let us denote by $O_K \subset [T_1,T_2]$ the set defined as $L_K$, but with the relation (\ref{condii}) replaced by
\begin{equation} \label{condiibis}
K^{-1} \nu^{-\beta} \leq -\min u_x \leq K \nu^{-\beta}.
\end{equation}

\begin{cor} \label{typicalcor}
For $K \geq K_1$ and $\nu < K_1^{-2/\beta}$, we have $\lambda(O_K) \geq C$.
\end{cor}

Now we fix
\begin{equation} \label{K}
K=K_1,
\end{equation}
and choose
\begin{equation} \label{nu0eq}
\nu_0=\Big( \frac{1}{6} K^{-2} \Big)^{1/\beta};\ C_1=\frac{1}{4}K^{-2};\ C_2=\frac{1}{20}K^{-4}.
\end{equation}
In particular, we have $0<C_1 \nu_0^{\beta} <C_2<1$: thus the intervals $J_i$ are non-empty and non-intersecting for all $\nu \in (0,\nu_0]$. Everywhere below the constants depend on $K$.
\\ \indent
Actually, we can choose any values of $C_1$, $C_2$ and $\nu_0$, provided:
\begin{equation} \label{nu0ineq}
C_1 \leq \frac{1}{4}K^{-2};\quad 5 K^2 \leq \frac{C_1}{C_2}<\frac{1}{\nu_0^{\beta}}.
\end{equation}

\begin{lemm} \label{lowerdiss}
For $\ell \in J_1$,
$$
S_{p}(\ell) \overset{p}{\gtrsim} \left\lbrace \begin{aligned} &\ell^{p},\ 0 \leq p \leq 1. \\ & \ell^{p} \nu^{-\beta(p-1)},\ p \geq 1. \end{aligned} \right.
$$
\end{lemm}

\begin{lemm} \label{lowerinert}
For $m \geq 0$ and $\ell \in J_2$,
$$
S_{p}(\ell) \overset{p}{\gtrsim} \left\lbrace \begin{aligned} & \ell^{p},\ 0 \leq p \leq 1. \\ & \ell,\ p \geq 1. \end{aligned} \right.
$$
\end{lemm}

Summing up the results above and the upper estimates proved in Section~\ref{turbupper} we obtain the following theorem.

\begin{theo} \label{avoir2}
For $\ell \in J_1$,
$$
S_{p}(\ell) \overset{p}{\sim} \left\lbrace \begin{aligned} & \ell^{p},\ 0 \leq p \leq 1. \\ & \ell^{p} \nu^{-(p-1) \beta},\ p \geq 1. \end{aligned} \right.
$$
On the other hand, for $\ell \in J_2$,
$$
S_{p}(\ell) \overset{p}{\sim} \left\lbrace \begin{aligned} & \ell^{p},\ 0 \leq p \leq 1. \\ & \ell,\ p \geq 1. \end{aligned} \right.
$$
\end{theo}

The following result follows immediately from the definition (\ref{flatness}).

\begin{cor} \label{flatnesscor}
For $\ell \in J_2$, the flatness satisfies $F(\ell) \sim \ell^{-1}$.
\end{cor}

\begin{theo} \label{spectrinert}
For $k$ such that $k^{-1} \in J_2$, we have $E(k) \sim k^{-2}$.
\end{theo}

\section{Acknowledgements}

I would like to thank P.~Clavin, A.~Kiselev, S.Kuksin, V.~Vicol and J.~Vovelle for helpful discussions. 

\section*{Appendix: the well-posedness of the fractional Burgers equation}

In this Appendix, we give a detailed proof for the well-posedness of the equation (\ref{fracBurgers}).  This proof is similar to the less detailed one given in \cite{KNS08}. We use the regularising effect of the Laplacian and the concept of mild solutions (i.e. Duhamel's formula). It can be generalised to the case of a stochastic forcing and/or a multidimensional setting: see \cite[Appendix A]{BorPhD} and \cite[Appendix 1]{BorM}. This proof cannot be adapted to the critical case $\alpha=1$, where more involved arguments using, for instance, a modulus of continuity, are needed: see \cite{CoVi12, KNS08}.
\\ \indent
Here, the functions whose Sobolev norms we consider do not necessarily have zero mean value in space. The only thing that changes is that now in the expressions for the Sobolev norms $W^{m,p}$ (resp. $H^s$) we have to add the norm in $L_p$ (resp. $L_2$) to the formulas in Section~\ref{Sob}. We use the standard notation $C(I,W^{m,p})$ for the space of continuous functions defined on the time interval $I$ with values in $W^{m,p}$ endowed with the supremum norm. The space $C(I,C^{\infty})$ will denote the intersection
$$
\cap_{m \geq 0}{C(I,H^m)}.
$$
\indent
We begin by considering mild solutions in $H^{1}$, in the spirit of \cite{DZ92,DZ96}. Then, by a bootstrap argument, we prove that for strictly positive times these solutions are actually smooth. Finally upper estimates (cf. Section~\ref{upper}) allow us to prove that such mild solutions are global.
\\ \indent
We will use the fact that the initial condition $u_0$ and the function $f$ in the nonlinearity are $C^{\infty}$-smooth.
\\ \indent
By a scaling argument, we can restrict ourselves to the equation (\ref{fracBurgers}) with $\nu=1$. We will denote by $S_L(t)$ the fractional heat semigroup $e^{-t \Lambda^{\alpha}}$. We recall that for $v \in L_2$ the function $S_L(t) v(\ix)$ is given by:
\begin{equation} \label{heat}
S_L(t) v(\ix)=\sum_{\ka \in \Z^d}{e^{-(2 \pi |\ka|)^{\alpha} t} \hat{v}_{\ka} e^{2 \pi i \ka \cdot \ix}}.
\end{equation}

We consider a mild form of (\ref{fracBurgers}):
\begin{equation} \label{mildHJ}
Y(t)=S_L(t) u_0+\int_{0}^{t}{\ S_L(t-\tau)\ (f(Y(\tau)))_x \ d \tau}.
\end{equation}
The fractional heat semigroup $S_L$ defines a contraction in each Sobolev space $H^s$. On the other hand, we have the following lemma.

\begin{lemm} \label{Lip}
The mapping
$$
Z \mapsto f(Z):\ H^1 \rightarrow H^1
$$
is locally Lipschitz on bounded subsets of $H^1$.
\end{lemm}

\textbf{Proof:} 
It suffices to develop $(f(Z_1)-f(Z_2))^{(1)}$ using Leibniz's formula and then to use the fact that by (\ref{Sobelement}), $\vert Z \vert_{\infty} \leq \Vert Z \Vert_1$.

\begin{lemm} \label{heatconvol}
For any $s \geq 0$, the operator
$$
Z \mapsto \Big( t \mapsto \int_{0}^{t}{S_L(t-\tau) Z(\tau) d \tau} \Big)
$$
maps bounded subsets of $C([0,T),H^{s})$ into bounded subsets of
\\
$C([0,T),H^{(s+(\alpha+1)/2)})$.
\end{lemm}

\textbf{Proof:}
Fix $s \geq 0$. By (\ref{heat}), for $\tau \in [0,t)$ we have
\begin{align} \nonumber 
& \Vert S_L(t-\tau) Z(\tau) \Vert^2_{s+(\alpha+1)/2}
\\ \nonumber
& \sim |(\hat{Z}(\tau))_\zero|^2+\sum_{\ka \in \Z^d}{|\ka|^{2s+\alpha+1} e^{-(2 \pi |\ka|)^{\alpha} (t-\tau) } |(\hat{Z}(\tau))_\ka|^2}
\\ \nonumber
& \lesssim |(\hat{Z}(\tau))_\zero|^2+\Big(\max_{\ka' \in \Z^d}{|\ka'|^{\alpha+1} e^{-(2 \pi |\ka'|)^{\alpha} (t-\tau) } } \Big) \sum_{\ka \in \Z^d}{|\ka|^{2s} |(\hat{Z}(\tau))_\ka|^2}
\\ \nonumber
& \lesssim \Big(1+\max_{\ka' \in \Z^d}{|\ka'|^{\alpha+1}e^{-(2 \pi |\ka'|)^{\alpha} (t-\tau) } } \Big) \Vert Z(\tau)\Vert^2_s.
\\ \nonumber 
& \lesssim C \Big[1+(t-\tau)^{-(\alpha+1)/\alpha} \Big] \Vert Z(\tau)\Vert^2_s.
\end{align}
To prove the lemma's statement, it remains to observe that since $\alpha>1$, we have:
$$
\int_{0}^{t}{(1+(t-\tau)^{-(\alpha+1)/\alpha})^{1/2} d \tau}<+\infty.
$$
\medskip
\\ \indent
Lemma~\ref{Lip}, Lemma~\ref{heatconvol} for $s=1$ and the Cauchy-Lipschitz theorem imply that the equation (\ref{mildHJ}) has a unique local solution in $H^{1}$.
\\ \indent
Now consider such a solution $Y$. We want to prove that this solution belongs to $C^{\infty}$ for all $t>0$. For this, since $\alpha>1$, it suffices to prove that for $s \geq 1$, a solution $Y \in H^s$ lies in the space $H^{s+(\alpha-1)/2}$. We will need the following result:

\begin{lemm} \label{Nemit}
For $s \geq 1$, the mapping
$$
Z \mapsto f(Z):\ H^{s} \rightarrow H^{s}
$$
is bounded on bounded subsets of $H^{s}$.
\end{lemm}

\textbf{Proof:} See \cite[Lemma A.0.5]{BorPhD}

\begin{theo} \label{bootstrap}
Consider a local solution $Y$ of (\ref{fracBurgers}) in $H^{1}$ defined on an interval $[0,T)$. If for some $s \geq 1$, we have $Y \in C([0,T),H^s)$, then we actually have $Y \in C([0,T),H^{s+(\alpha-1)/2})$.
\end{theo}

\textbf{Proof:} By Lemma~\ref{Nemit} we have
$$
(f(Y(\tau)))_x \in C([0,T),H^{s-1}),
$$
and thus by Lemma~\ref{heatconvol} we get
$$
\int_{0}^{t}{S_L(t-\tau)(f(Y(\tau)))_x  d\tau} \in C([0,T),H^{(s+(\alpha-1)/2)}).
$$
Since $Y$ is a solution of (\ref{mildHJ}) and the semigroup $S_L$ is smoothing,
$$
Y(t)=S_L(t) u_0+\int_{0}^{t}{S_L(t-\tau) (f(Y(\tau)))_x  d \tau}
$$
belongs to the space $C([0,T),H^{s+(\alpha-1)/2}).$
\medskip
\\ \indent
Thus, by a bootstrap argument, it follows that there exists a unique local solution to (\ref{fracBurgers}), which is $C^{\infty}$-smooth in space for $t>0$. To prove that this solution is necessarily global, it suffices to use the uniform in time upper estimates in Section~\ref{upper}.\\ \indent

\bibliographystyle{plain}
\bibliography{Bibliogeneral}

\end{document}